\documentclass[12pt]{amsart}
\usepackage{fullpage}
\usepackage{enumerate}

\usepackage{fullpage}

\usepackage{amsmath,amsthm,amsfonts,amssymb,amscd,amsopn,amsmath,mathrsfs,latexsym}

\usepackage[dvipsnames]{xcolor}

\usepackage[all]{xy}
\usepackage[colorlinks]{hyperref}

\newdimen{\PDww}
\newdimen{\PDhh}
\newdimen{\PDtt}
\newcounter{PDcounter}%
\def\PD{\futurelet\next\PDfindnfoo}
\def\PDfindnfoo{\ifx\next[%
    \expandafter\PDfindn%
\else%
    \PDfindn[1]%
\fi}
\def\PDfindn[#1]{\def\PDn{#1}\expandafter\PDfindnfoofoo}
\def\PDfindnfoofoo{\futurelet\next\PDfindwfoo}

\def\PDfindwfoo{\ifx\next[%
    \expandafter\PDfindw%
\else%
    \PDfindw[1pc]%
\fi}
\def\PDfindw[#1]{\setlength{\PDww}{#1}\expandafter\PDfindwfoofoo}
\def\PDfindwfoofoo{\futurelet\next\PDfindhfoo}

\def\PDfindhfoo{\ifx\next[%
    \expandafter\PDfindh%
\else%
    \PDfindh[1.25pc]%
\fi}
\def\PDfindh[#1]{\setlength{\PDhh}{#1}\expandafter\PDfindhfoofoo}
\def\PDfindhfoofoo{\futurelet\next\PDfindtfoo}

\def\PDfindtfoo{\ifx\next[%
    \expandafter\PDfindt%
\else%
    \PDfindt[0.2pt]%
\fi}
\def\PDfindt[#1]{\setlength{\PDtt}{#1}\expandafter\PDcontent}
\long\def\PDcontent#1{%
\def\select##1##2##3{%
    \let\xtxt##1\let\ytxt##2%
    \if\xtxt\ytxt##3\fi}
\def\selectcolor##1{%
    \color{black}%
    \select{##1}{0}{\color{gray}}%
    \select{##1}{1}{\color{red}}%
    \select{##1}{2}{\color{blue}}%
    \select{##1}{3}{\color{Green}}%
    \select{##1}{4}{\color{orange}}%
    \select{##1}{5}{\color{cyan}}%
    \select{##1}{6}{\color{magenta}}%
    \select{##1}{7}{\color{teal}}%
    \select{##1}{8}{\color{brown}}%
    \select{##1}{9}{\color{violet}}%
    }
\def\PDvrule##1##2##3{
    \rule{-0.5\dimexpr##2}{0pc}%
    \rule[\dimexpr##1]{##2}{##3}%
    \rule{-0.5\dimexpr##2}{0pc}%
    }%
\def\PDhhrule##1##2##3{
    \rule[\dimexpr##1-0.5\dimexpr##3]{##2}{##3}%
    }%
\def\PDframe{%
\raisebox{0pc}[\PDhh][0pc]{\makebox[0pc][l]{%
    \PDvrule{0pc}{4\PDtt}{\PDhh}%
    {\color{lightgray}%
        \setcounter{PDcounter}{0}%
        \loop\ifnum\value{PDcounter}<\PDn\relax%
            \hspace{\PDww}%
            \PDvrule{0pc}{0.25\PDtt}{\PDhh}%
            \addtocounter{PDcounter}{1}%
        \repeat}%
    \PDvrule{0pc}{4\PDtt}{\PDhh}%
    \rule{-\PDn\PDww}{0pt}%
    \PDhhrule{\PDhh}{\PDn\PDww}{4\PDtt}%
    \rule{-\PDn\PDww}{0pt}%
    \PDhhrule{0pc}{\PDn\PDww}{4\PDtt}%
    }}}
\def\T##1{%
    {\let\F\FinT%
    \let\J\JinT%
    \let\L\LinT%
    \let\G\GinT%
    \let\I\IinT%
    \let\H\HinT%
    \let\O\OinT%
    \PDframe%
    \makebox[\PDn\PDww][l]{##1}}%
}
\def\FinT##1{%
    \makebox[0pc][l]{\selectcolor{##1}%
    \rule{\dimexpr##1\PDww-0.5\PDww}{0pc}%
    \PDvrule{0pc}{6\PDtt}{0.5\PDhh}%
    \PDhhrule{0.5\PDhh}{\dimexpr\PDn\PDww-##1\PDww+0.5\PDww}{6\PDtt}%
    }}
\def\JinT##1{%
    \makebox[0pc][l]{\selectcolor{##1}%
    \PDhhrule{0.5\PDhh}{\dimexpr##1\PDww-0.5\PDww}{6\PDtt}%
    \PDvrule{0.5\PDhh}{6\PDtt}{0.5\PDhh}%
    }}
\def\IinT##1{%
    \makebox[0pc][l]{\selectcolor{##1}%
    \rule{\dimexpr##1\PDww-0.5\PDww}{0pc}%
    \PDvrule{0pc}{6\PDtt}{\PDhh}%
    }}
\def\HinT##1{%
    \makebox[0pc][l]{\selectcolor{##1}%
    \PDhhrule{0.5\PDhh}{\dimexpr\PDn\PDww}{6\PDtt}%
    }}
\def\OinT##1##2{%
    \makebox[0pc][l]{\selectcolor{##1}%
    \rule{\dimexpr##1\PDww-0.5\PDww}{0pc}%
    \raisebox{0.25\PDhh}{\makebox[0pc][c]{\(\displaystyle##2\)}}%
    }}
\def\GinT##1{%
    \makebox[0pc][l]{\selectcolor{##1}%
    \PDhhrule{0.5\PDhh}{\dimexpr##1\PDww-0.5\PDww}{6\PDtt}%
    \PDvrule{0pc}{6\PDtt}{0.5\PDhh}%
    }}
\def\LinT##1{%
    \makebox[0pc][l]{\selectcolor{##1}%
    \rule{\dimexpr##1\PDww-0.5\PDww}{0pc}%
    \PDvrule{0.5\PDhh}{6\PDtt}{0.5\PDhh}%
    \PDhhrule{0.5\PDhh}{\dimexpr\PDn\PDww-##1\PDww+0.5\PDww}{6\PDtt}%
    }}
\def\x{%
    \raisebox{\dimexpr0.5\PDhh-0.25pc}{\makebox[0pc][c]{\color{gray}\(\displaystyle\bullet\)}}%
    }
\def\o{%
    \raisebox{\dimexpr0.5\PDhh-0.25pc}{\makebox[0pc][c]{\color{gray}\(\displaystyle\circ\)}}%
    }
\def\MM##1{%
\raisebox{0.25\PDhh}[\PDhh]{%
    \def\e{}\def\inputtxt{##1}%
    \ifx\e\inputtxt\relax\else%
        \setcounter{PDcounter}{1}%
        \makebox[\PDn\PDww][l]{%
            \expandafter\Mbegin##1\Mend%
        }%
    \fi}}
\def\Mbegin##1{%
    \ifx\Mend##1\relax\else%
        {\selectcolor{\thePDcounter}%
            \makebox[\PDww]{\(\displaystyle ##1\)}}%
        \addtocounter{PDcounter}{1}
        \expandafter\Mbegin
    \fi}
\def\subframe{%
    \raisebox{0pc}[\PDhh][0pc]{\makebox[0pc][l]{%
    \color{lightgray}%
    \PDvrule{0pc}{\PDtt}{\PDhh}%
    \hspace{\PDww}%
    \PDvrule{0pc}{\PDtt}{\PDhh}%
    \rule{-\PDww}{0pt}%
    \PDhhrule{\PDhh}{\PDww}{\PDtt}%
    \rule{-\PDww}{0pt}%
    \PDhhrule{0pc}{\PDww}{\PDtt}%
    }}}
\def\F##1{\subframe\makebox[\PDww]{\selectcolor{##1}%
    \PDvrule{0pc}{6\PDtt}{0.5\PDhh}%
    \PDhhrule{0.5\PDhh}{0.5\PDww}{6\PDtt}%
    \rule{-0.5\PDww}{0pt}%
    }}%
\def\J##1{\subframe\makebox[\PDww]{\selectcolor{##1}%
    \PDvrule{0.5\PDhh}{6\PDtt}{0.5\PDhh}%
    \rule{-0.5\PDww}{0pt}%
    \PDhhrule{0.5\PDhh}{0.5\PDww}{6\PDtt}%
    }}%
\def\L##1{\subframe\makebox[\PDww]{\selectcolor{##1}%
    \PDvrule{0.5\PDhh}{6\PDtt}{0.5\PDhh}%
    \PDhhrule{0.5\PDhh}{0.5\PDww}{6\PDtt}%
    \rule{-0.5\PDww}{0pt}%
    }}%
\def\G##1{\subframe\makebox[\PDww]{\selectcolor{##1}%
    \PDvrule{0\PDhh}{6\PDtt}{0.5\PDhh}%
    \rule{-0.5\PDww}{0pt}%
    \PDhhrule{0.5\PDhh}{0.5\PDww}{6\PDtt}%
    }}%
\def\H##1{\subframe\makebox[\PDww]{\selectcolor{##1}%
    \PDhhrule{0.5\PDhh}{\PDww}{6\PDtt}%
    }}%
\def\I##1{\subframe\makebox[\PDww]{\selectcolor{##1}%
    \PDvrule{0pc}{6\PDtt}{\PDhh}%
    }}%
\def\X##1##2{\subframe\makebox[\PDww]{\selectcolor{##1}%
    \rule{-0.5\PDww}{0pt}%
    \PDhhrule{0.5\PDhh}{\PDww}{6\PDtt}%
    \selectcolor{##2}%
    \rule{-0.5\PDww}{0pt}%
    \PDvrule{0pc}{6\PDtt}{\PDhh}%
    }}%
\def\O##1##2{\subframe\makebox[\PDww]{\selectcolor{##1}%
    \raisebox{0.25\PDhh}{\makebox[0pc][c]{\(\displaystyle##2\)}}
    }}%
\def\M##1##2{\makebox[\PDww]{\selectcolor{##1}%
    \raisebox{0.25\PDhh}[\PDhh]{\makebox[0pc][c]{\(\displaystyle##2\)}}
    }}%
\def\,{\rule{\PDww}{0pc}\rule{0pc}{\PDhh}}%
{\def\arraystretch{0}\setlength{\arraycolsep}{0pc}
\begin{array}{l}
    #1
\end{array}}}

\def\HPD{\operatorname{HPD}}
\def\wt{\operatorname{wt}}

\newtheorem{Th}{Theorem}[section]

\newtheorem{Coro}[Th]{Corollary}

\newtheorem{Eg}[Th]{Example}

\author{Yihan Xiao}
\address[Yihan Xiao]{School of science, Shanghai university, 99 Shangda Rd, Bao Shan Qu, Shang Hai Shi, China, 200444}
\email{jstbfrnd@shu.edu.cn}

\author{Rui Xiong}
\address[Rui Xiong]{Department of Mathematics and Statistics, University of Ottawa, 150 Louis-Pasteur, Ottawa, ON, K1N 6N5, Canada}
\email{rxion043@uottawa.ca}

\author{Haofeng Zhang}
\address[Haofeng Zhang]{School of mathematical science, Peking university, 5 Yiheyuan Rd, Hai Dian Qu, Bei Jing Shi, China, 100871}
\email{aprilgrimoire@stu.pku.edu.cn}

\title{Hybrid Pipe Dreams for Key Polynomials}

\date{\today}

\begin{document}

\maketitle

\begin{abstract}
We develop a family of new combinatorial models for key polynomials.
It is similar to the hybrid pipe dream model for Schubert polynomials defined recently by Knutson and Udell. 
\end{abstract}

\section{Introduction}

The key polynomials $\kappa_\alpha(x)$ for $\alpha\in \mathbb{Z}_{\geq0}^n$ are the Demazure characters for general linear groups $GL_n$. 
They are a non-symmetric generalization of the Schur polynomials $s_\lambda(x_1,\ldots,x_n)$, the Weyl characters for $GL_n$.
Thanks to Demazure \cite{Dem1974A,Dem1974B}, the key polynomials are characterized via Demazure operators $\pi_i$ ($1\leq i\leq n$) on the polynomial ring $\mathbb{Z}[x_1,\ldots,x_n]$ by 
$$\pi_i f = \frac{x_i\cdot f-x_{i+1}\cdot f|_{x_i\leftrightarrow x_{i+1}}}{x_i-x_{i+1}}.$$
Precisely, when $\alpha=(\alpha_1,\ldots,\alpha_n)$ is a partition, i.e., $\alpha_1\geq \cdots\geq \alpha_n$, the key polynomial is a monomial
$$\kappa_\alpha(x)=x_1^{\alpha_1}\cdots x_n^{\alpha_n}.$$
Otherwise, if $\alpha_i<\alpha_{i+1}$, we have 
$$\kappa_{\alpha} = \pi_i\kappa_{\alpha'}$$
where $\alpha'$ is the composition obtained from $\alpha$ by exchanging $\alpha_i$ and $\alpha_{i+1}$. 

The characterization above determines key polynomials algebraically, while combinatorial formulas are needed for computation and study their properties. 
The first combinatorial object for key polynomials is due to Lascoux and Sch\"utzenberger \cite{LS90} by a generalization of semi-standard Young tableaux.
After that, many different combinatorial models were introduced 
\cite{Ko91,RS95,As18,AS18,AQ19,BS20} by 
Kohnert, Reiner, Shimozono, Assaf, Quijad Searles, Busiumas, Scrimshaw, etc. 
For general Lie types, combinatorial models of Demazure characters were studied in
\cite{Li94,Li95,Ka93,LP07,LP08} by Littlemann, Kashiwara, Lenart, Postinkov, etc.

On the other hand, Macdonald \cite{Mac88}
introduced a family of symmetric functions $P_\lambda$ with two variables $q$, $t$, see Macdonald's book \cite[chapter VI]{Mac95}.  
By (4.14) loc. cit., they specialize to Schur functions by setting $q = t$. 
The non-symmetric generalization $E_\alpha(x;q,t)$ was developed
by Opdam \cite{Op95}, Macdonald \cite{Mac96} and Cherednik \cite{Ch95}.
It was first pointed out by Ion \cite{Io03} that 
$$\kappa_\alpha(x)=E_\alpha(x;q=\infty,t=\infty).$$
Haglund, Haiman and Loehr \cite{HHL08} gave the first combinatorial model of $E_\alpha$ in terms of non-attacking skyline fillings. 
In particular, it gives a new combinatorial model for key polynomials, studied by \cite{Mo17,MPS21,As18} by Monical, Monical, Pechenik, Searles, Assaf, etc.

There are many other combinatorial models of $E_\alpha$, see \cite{RY11,BW22,ABW23,Le08,Al19} by Ram, Yip, Borodin, Wheeler, Aggarwal, Lenart, Alexandersson, etc. 
Each of them should specialize a combinatorial model of key polynomials. 
The model we are going to specialize in this paper is from Borodin and Wheeler \cite{BW22} via vertex models. 
We remark that in the proof of loc. cit., we have to keep two variables $q,t$ unspecialized, so we do not have a good explanation why this model satisfies the algebraic recursion given by Demazure operators.  
Besides, the result of Busiumas and Scrimshaw \cite{BS20} implies another vertex model for key polynomials, which agrees with the one specialized from \cite{ABW23}.

Schur polynomials admit another family of non-symmetric generalization known as Schubert polynomials $\mathfrak{S}_w$ for permutations $w$. 
Schubert polynomials were also introduced by Lascoux and Sch\"utzenberger \cite{LS82} to study the intersection theory over flag varieties. 
A modern treatment of this topic can be found in the book by Anderson and Fulton \cite{AF}. 
They satisfy a similar recursion formula as key polynomials, using a different family of operators. 
Based on Fomin and Stanley \cite{FS94} and Fomin and Kirillov \cite{FK96}, Bergeron and Billey \cite{BS93} defined a combinatorial model for Schubert polynomials, now known as pipe dreams (PD). 
Many years after PD was found, Lam, Lee and Shimozono \cite{LLS21} discovered a different combinatorial model for Schubert polynomials called bumpless pipe dreams (BPD).
Although the two models compute the same polynomial, it is far from obvious that there exists a weight-preserving bijection between them. 
A ``canonical'' bijection was established by Gao and Huang \cite{GH}. 

Recently, a new approach was introduced by Knutson and Udell \cite{KU24}. 
They gave a new family of models called hybrid pipe dreams (HPD) for Schubert polynomials including both PD and BPD.
In their model, any two of HPD models are connected by a path of HPD models with explicit local weight-preserving bijections at each step. 
In particular, it gives a bijection between PD and BPD inductively. 

In this paper, we develop the hybrid version for key polynomials (Theorem \ref{thm:mainthm}). 
It generalizes the model specialized from \cite{BW22}. 
As an application, we derive a two-side branching rule for key polynomials (Corollary \ref{coro:branch}). 

Let us briefly discuss our approach. 
First, we take the model specialized from Borodin and Wheeler \cite{BW22} as an initial model, which should be viewed as an analogy of PD for Schubert polynomials. 
Then, by introducing another family of pieces, we can define our hybrid pipe dream model. 
Lastly, we find an explicit weight-preserving bijection between different models, similar to Knutson and Udell \cite[\textsection 2]{KU24}. 
The most essential step in the proof is the case of exchanging two adjacent rows of different types.

Let us finish the introduction by discuss some further questions. 
It is natural to ask whether our model can be extended to non-symmetric Macdonald polynomials $E_\alpha$. 
But the weights for $E_\alpha$ in Borodin and Wheeler \cite{BW22} are quite complicated, and we do not have a good guess for its hybrid version. 

It was noticed by Lascoux and Sch\"utzenberger \cite{LS90B} that 
Schubert polynomials expand positively into key polynomials. 
A proof can be found in Reiner and Shimozono \cite{RS95}, 
see also Assaf and Schilling \cite{AS18B}, Shimozono and Yu \cite{SY23} and Gold, Mili\'cevi\'c and Sun \cite{GMS24}. 
It is natural to ask whether such relations can be seen from HPD models for Schubert polynomials and key polynomials.

Recently, Yu \cite[Theorem 6.7]{Yu24} gave a bijection between HPD for Schubert polynomials to Bruhat chains. 
For key polynomials, a chain formula is also expected. Our model should give some clue in this direction. 

There are many conjectures in \cite{MTY19} concerning key polynomials. 
Some of them were confirmed recently in \cite{FMD18,FG21,FGPS20} by Fink, M\'esz\'aros, Dizier, Fan, Guo, Peng, Sun, etc. 
Since our model has the feature that the rows contributing $x_i$ and $x_j$ can be adjacent for any not necessarily adjacent indices $i\neq j$, it possibly provides a way of approaching these conjectures.

The paper is organized as follows. 
In section \ref{subsec:model} and \ref{subsec:mainthm} we introduce our model and state the main theorem (Theorem \ref{thm:mainthm}). 
In section \ref{subsec:branch} we discuss the two-side branching rule (Corollary \ref{coro:branch}). 
In section \ref{subsec:Preproof}, we identify the model specialized from Borodin and Wheeler \cite{BW22} with one of our model, and reduce the proof of Theorem \ref{thm:mainthm} to the case of swapping two adjacent rows. 
Section \ref{sec:swaprow} will be devoted to establish the swapping. 
For this purpose, we need to distinguish two types of sub-columns, and they are dealt separately in section \ref{subsec:frozen} and \ref{subsec:unfrozen}. 

\subsection*{Acknowledgment}
Part of this work was carried out during the PACE program in the summer of 2024 at BICMR in Peking University. 
We would also like to
thank Neil J.Y. Fan, Yibo Gao, Peter L. Guo, Mark Shimozono and Tianyi Yu for helpful discussions.

\section{Pipe dream model}\label{sec:PDmod}


\subsection{The tiles}
\label{subsec:model}
In this paper, we will fix an integer $n\geq 1$.
we will denote $[n]=\{1,\ldots,n\}$, and take it as the set of colors.
Our pipe dream model consists of two types of tiles which can be decomposed further into smaller sub-tiles. 

A \emph{west tile} is a horizontal stack of $n$ sub-tiles 
$$
\PD[1][1pc][1.2pc]{\J0}\qquad 
\PD[1][1pc][1.2pc]{\F0}\qquad 
\PD[1][1pc][1.2pc]{\H0}\qquad 
\PD[1][1pc][1.2pc]{\I0}\qquad 
\PD[1][1pc][1.2pc]{\X00}\qquad 
\PD[1][1pc][1.2pc]{\O0{}}$$
with each ``pipe'' colored by a number from $[n]$. 
We require the following conditions on the west tiles
\begin{itemize}
    \item only the pipe colored by $i\in [n]$ can touch the $i$-th floor or the $i$-th ceiling;
    \item 
    in the sub-tile $\PD{\X00}$, the horizontal pipe is colored less than the vertical pipe. 
\end{itemize}
In particular, a pipe can only turn at most once in a tile.

Similarly, an \emph{east tile} is a horizontal stack of $n$ sub-tiles 
$$
\PD[1][1pc][1.2pc]{\L0}\qquad 
\PD[1][1pc][1.2pc]{\G0}\qquad 
\PD[1][1pc][1.2pc]{\H0}\qquad 
\PD[1][1pc][1.2pc]{\I0}\qquad 
\PD[1][1pc][1.2pc]{\X00}\qquad 
\PD[1][1pc][1.2pc]{\O0{}}$$
with each ``pipe'' colored by a number from $[n]$. 
We require the following conditions on the east tiles
\begin{itemize}
    \item only the pipe colored by $i\in [n]$ can touch the $i$-th floor or the $i$-th ceiling;
    \item 
    in the sub-tile $\PD{\X00}$, the horizontal pipe is colored larger than the vertical pipe. 
\end{itemize}
Note that conditions on the crossing of pipes are different.

\begin{Eg}
When $n=5$, here are some examples of west tiles. 
$$\PD[5]{
\,\MM{{}2{}45}\,\\
\,\T{\I2\I4\I5}\\
\,\MM{{}2{}45}}
\quad 
\PD[5]{
\,\MM{{}2{}{}5}\,\\
\M11\T{\H1\I2\I5}\M11\\
\,\MM{{}2{}{}5}}
\quad 
\PD[5]{
\,\MM{1{}{}{}5}\,\\
\,\T{\F3\I1\I5}\M3{3}\\
\,\MM{1{}3{}5}\,}
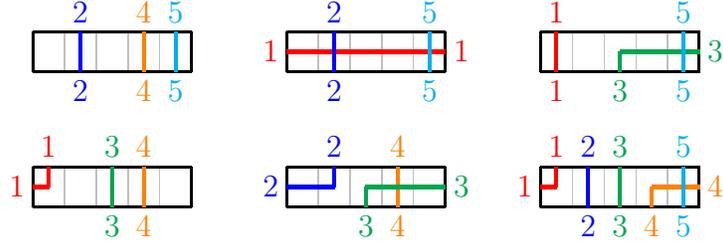
$$
$$
\PD[5]{
\,\MM{1{}34{}}\,\\
\M1{1}\T{\J1\I3\I4}\\
\,\MM{{}{}34}}
\quad 
\PD[5]{
\,\MM{{}2{}4{}}\,\\
\M2{2}\T{\J2\I4\F3}\M3{3}\\
\,\MM{{}{}34}}
\quad 
\PD[5]{
\,\MM{123{}5}\,\\
\M1{1}\T{\J1\I2\F4\I3\I5}\M4{4}\\
\,\MM{{}2345}}
$$    
Here are some examples of east tiles. 
$$\PD[5]{
\,\MM{{}2{}45}\,\\
\,\T{\I2\I4\I5}\\
\,\MM{{}2{}45}}
\quad 
\PD[5]{
\,\MM{1{}3}\,\\
\M44\T{\H4\I3\I1}\M44\\
\,\MM{1{}3}}
\quad 
\PD[5]{
\,\MM{1{}{}45}\,\\
\M22\T{\G2\I1\I4\L5}\M55\\
\,\MM{12{}4}}
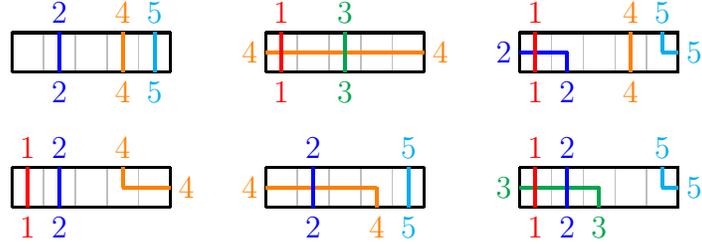
$$
$$
\PD[5]{
\,\MM{12{}4}\,\\
\,\T{\L4\I2\I1}\M44\\
\,\MM{12}\,}
\quad 
\PD[5]{
\,\MM{{}2{}{}5}\,\\
\M44\T{\G4\I5\I2}\,\\
\,\MM{{}2{}4{5}}}
\quad 
\PD[5]{
\,\MM{12{}{}5}\,\\
\M33\T{\L5\G3\I1\I2}\M55\\
\,\MM{123}}
$$
\end{Eg}

\subsection{Main theorem}
\label{subsec:mainthm}
Let $\alpha=(\alpha_1,\ldots,\alpha_n)$ be a composition.
We will fix an integer $N$ such that $\alpha_i\leq N$ for all $i$. 
We will make a choice of type 
$\tau=(\tau_1,\ldots,\tau_n)$ with $\tau_i\in \{W,E\}$.

Consider a grid of $n\times (N+1)$ as follows:
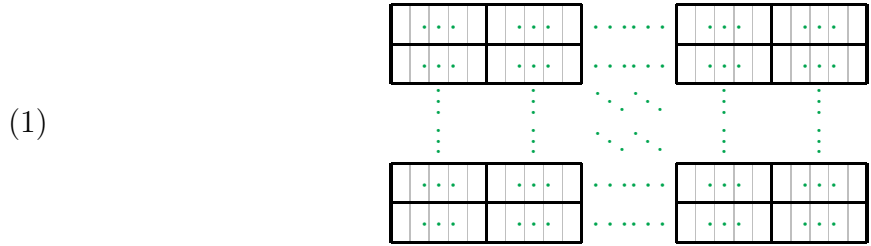
\begin{equation}\label{eq:grid}
\PD[5][0.6pc][1.25pc]{
\T{\O3{\cdots}}\T{\O3{\cdots}}
\M0{}\M3{\cdots}\M0{}\M3{\cdots}\M0{}
\T{\O3{\cdots}}\T{\O3{\cdots}}\\
\T{\O3{\cdots}}\T{\O3{\cdots}}
\M0{}\M3{\cdots}\M0{}\M3{\cdots}\M0{}
\T{\O3{\cdots}}\T{\O3{\cdots}}\\
\MM{{}{}{\vdots}}\MM{{}{}{\vdots}}
\M0{}\M3{\ddots}\M0{}\M3{\ddots}\M0{}
\MM{{}{}{\vdots}}\MM{{}{}{\vdots}}\\
\MM{{}{}{\vdots}}\MM{{}{}{\vdots}}
\M0{}\M3{\ddots}\M0{}\M3{\ddots}\M0{}
\MM{{}{}{\vdots}}\MM{{}{}{\vdots}}\\
\T{\O3{\cdots}}\T{\O3{\cdots}}
\M0{}\M3{\cdots}\M0{}\M3{\cdots}\M0{}
\T{\O3{\cdots}}\T{\O3{\cdots}}\\
\T{\O3{\cdots}}\T{\O3{\cdots}}
\M0{}\M3{\cdots}\M0{}\M3{\cdots}\M0{}
\T{\O3{\cdots}}\T{\O3{\cdots}}
}
\end{equation}
We refer to the leftmost column as the $0$-th column, and the next one as the $1$-st column, etc. So the rightmost column is the $N$-th column. 
We refer to the top row as the first row, and the next one as the second row, etc. So the bottom row is the $n$-th row. 
We say the $i$-th row is of west (resp., east) type, if $\tau_i$ is $W$ (resp., $E$). 

Next, we label the colors $[n]$ on the boundary of the grid. 
Let $k$ be the number of rows of west type. 
On the left (resp., right) boundary, we label those rows of west (resp., east) type by colors $1,2,\cdots,k$ (resp., $n,n-1,\cdots,k+1$), from up to down. 
On the upper boundary, we label the color $i\in [n]$ on the $i$-th ceiling of the $\alpha_i$-th column. 
We left the lower boundary empty. 

We refer the row labeled by $i$ above by the $i$-row.
Namely, if $i\leq k$ (resp., $i>k$), then the $i$-row is the $i$-th (resp., $(n+1-i)$-th) row of west (resp., east) type. 
Note that the $i$-row is NOT the $i$-th row.

Let us denote by
$\HPD_{\tau}(\alpha)$
the set of tilings of the grid, with the rows of west (resp. east) type tiled by west (resp., east) tiles.
For each $\pi\in \HPD_{\tau}(\alpha)$, we define its weight 
\begin{equation}\label{eq:weight}
\wt(\pi)=x_1^{m_1}\cdots x_n^{m_n},
\end{equation}
where $m_i$ is the number of tiles in the $i$-row with a pipe touching (resp., not touching) the right wall if the $i$-row is of west (resp., east) type. 

The main theorem of this paper is the following.

\begin{Th}\label{thm:mainthm}
For any choice of $\tau\in \{W,E\}^n$, we have
$$\kappa_\alpha(x)=\sum_{\pi\in \HPD_{\tau}(\alpha)}\wt(\pi).$$
\end{Th}

\begin{Eg}For example, if $n=4$ and $\alpha=(1,3,0,2)$. 
we can take $N=3$ and $\tau=WEEW$. 
Consider the following element in $\HPD_{\tau}(\alpha)$
(we draw a dot on the walls contributing weights)
$$\PD[4]{
\,\,\,
\MM{{}{}{3}{}}\MM{{1}{}{}{}}\MM{{}{}{}{4}}\MM{{}{2}{}{}}\\
\M0{x_1^3}\,
\M11\T{\H1\I3}\x\T{\J1\F2}\x\T{\H2\I4}\x\T{\J2}\,\,
    \M0{W}\\
\M0{x_4^0}\,
\,\T{\L3}\T{\I2\H3}\T{\G3\L4}\T{\H4}\M44\,
    \M0{E}\\
\M0{x_3^1}\,
\,\T{}\x\T{\L2}\T{\G2\L3}\T{\H3}\M33\,
    \M0{E}\\
\M0{x_2^2}\,
\M22\T{\H2}\x\T{\H2}\x\T{\J2}\T{}\,\,
    \M0{W}}$$
Its weight is $x_1^3x_2^2x_3$. 
The following table might be helpful: 
$$\begin{array}{ccccl}\hline
\text{$i$-row} & \text{descrption} &
\text{type} &m_i&\hfill
\text{explanation} \hfill\\\hline
\text{$1$-row} & \text{the $1$-st row}
&W & 3
& \text{since pipe $1$ intersects once and pipe $2$ intersects twice}\\
\text{$2$-row} & \text{the $4$-th row}
& W & 2 & \text{since pipe $2$ intersects twice}\\
\text{$3$-row} & \text{the $3$-rd row}
& E & 1 & \text{since there is one wall without pipe intersected}\\
\text{$4$-row} & \text{the $2$-nd row}
& E & 0 & \text{since all walls are intersected}\\\hline
\end{array}$$
Besides this, there rest $12$ elements in $\HPD_{\tau}(\alpha)$, see 
Figure \ref{fig:1302}. 
The key polynomial is exactly the sum of their weights. 

\begin{figure}[hhhh]
$$\begin{matrix}
\PD[4][0.5pc]{
\T{\H1\I3}\x\T{\J1}\T{\I4}\T{\I2}\\
\T{\I3}\x\T{}\x\T{\L4}\T{\I2\H4}\\
\T{\L3}\T{\H3}\T{\H3}\T{\I2\H3}\\
\T{\H2}\x\T{\H2}\x\T{\H2}\x\T{\J2}}\\\vphantom{\dfrac{1}{2}}
x_1x_2^3x_4^2
\end{matrix}\quad 
\begin{matrix}
\PD[4][0.5pc]{
\T{\H1\I3}\x\T{\J1}\T{\I4}\T{\I2}\\
\T{\L3}\T{\G3}\x\T{\L4}\T{\I2\H4}\\
\T{}\x\T{\L3}\T{\H3}\T{\I2\H3}\\
\T{\H2}\x\T{\H2}\x\T{\H2}\x\T{\J2}}
\\\vphantom{\dfrac{1}{2}}
x_1x_2^3x_3x_4
\end{matrix}\quad 
\begin{matrix}
\PD[4][0.5pc]{
\T{\H1\I3}\x\T{\J1}\T{\I4}\T{\I2}\\
\T{\L3}\T{\H3}\T{\G3\L4}\T{\I2\H4}\\
\T{}\x\T{}\x\T{\L3}\T{\I2\H3}\\
\T{\H2}\x\T{\H2}\x\T{\H2}\x\T{\J2}}
\\\vphantom{\dfrac{1}{2}}
x_1x_2^3x_3^2
\end{matrix}\quad 
\begin{matrix}
\PD[4][0.5pc]{
\T{\H1\I3}\x\T{\J1\F4}\x\T{\J4}\T{\I2}\\
\T{\I3}\x\T{\L4}\T{\H4}\T{\I2\H4}\\
\T{\L3}\T{\H3}\T{\H3}\T{\I2\H3}\\
\T{\H2}\x\T{\H2}\x\T{\H2}\x\T{\J2}}
\\\vphantom{\dfrac{1}{2}}
x_1^2x_2^3x_4
\end{matrix}$$
$$ 
\begin{matrix}
\PD[4][0.5pc]{
\T{\H1\I3}\x\T{\J1\F4}\x\T{\J4}\T{\I2}\\
\T{\L3}\T{\G3\L4}\T{\H4}\T{\I2\H4}\\
\T{}\x\T{\L3}\T{\H3}\T{\I2\H3}\\
\T{\H2}\x\T{\H2}\x\T{\H2}\x\T{\J2}}
\\\vphantom{\dfrac{1}{2}}
x_1^2x_2^3x_3
\end{matrix}\quad
\begin{matrix}
\PD[4][0.5pc]{
\T{\H1\I3}\x\T{\J1}\T{\F2\I4}\x\T{\J2}\\
\T{\I3}\x\T{}\x\T{\I2\L4}\T{\H4}\\
\T{\L3}\T{\H3}\T{\I2\H3}\T{\H3}\\
\T{\H2}\x\T{\H2}\x\T{\J2}\T{}}
\\\vphantom{\dfrac{1}{2}}
x_1^2x_2^2x_4^2
\end{matrix}\quad
\begin{matrix}
\PD[4][0.5pc]{
\T{\H1\I3}\x\T{\J1}\T{\F2\I4}\x\T{\J2}\\
\T{\L3}\T{\G3}\x\T{\I2\L4}\T{\H4}\\
\T{}\x\T{\L3}\T{\I2\H3}\T{\H3}\\
\T{\H2}\x\T{\H2}\x\T{\J2}\T{}}
\\\vphantom{\dfrac{1}{2}}
x_1^2x_2^2x_3x_4
\end{matrix}\quad
\begin{matrix}
\PD[4][0.5pc]{
\T{\H1\I3}\x\T{\J1}\T{\F2\I4}\x\T{\J2}\\
\T{\L3}\T{\H3}\T{\I2\G3\L4}\T{\H4}\\
\T{}\x\T{}\x\T{\I2\L3}\T{\H3}\\
\T{\H2}\x\T{\H2}\x\T{\J2}\T{}}
\\\vphantom{\dfrac{1}{2}}
x_1^2x_2^2x_3^2
\end{matrix}$$
$$
\begin{matrix}
\PD[4][0.5pc]{
\T{\H1\I3}\x\T{\J1\F2}\x\T{\H2\I4}\x\T{\J2}\\
\T{\I3}\x\T{\I2}\x\T{\L4}\T{\H4}\\
\T{\L3}\T{\I2\H3}\T{\H3}\T{\H3}\\
\T{\H2}\x\T{\J2}\T{}\T{}}
\\\vphantom{\dfrac{1}{2}}
x_1^3x_2x_4^2
\end{matrix}\quad
\begin{matrix}
\PD[4][0.5pc]{
\T{\H1\I3}\x\T{\J1\F2}\x\T{\H2\I4}\x\T{\J2}\\
\T{\I3}\x\T{\L2}\T{\G2\L4}\T{\H4}\\
\T{\L3}\T{\H3}\T{\I2\H3}\T{\H3}\\
\T{\H2}\x\T{\H2}\x\T{\J2}\T{}}
\\\vphantom{\dfrac{1}{2}}
x_1^3x_2^2x_4
\end{matrix}\quad 
\begin{matrix}
\PD[4][0.5pc]{
\T{\H1\I3}\x\T{\J1\F2}\x\T{\H2\I4}\x\T{\J2}\\
\T{\L3}\T{\I2\G3}\x\T{\L4}\T{\H4}\\
\T{}\x\T{\I2\L3}\T{\H3}\T{\H3}\\
\T{\H2}\x\T{\J2}\T{}\T{}}
\\\vphantom{\dfrac{1}{2}}
x_1^3x_2x_3x_4
\end{matrix}\quad 
\begin{matrix}
\PD[4][0.5pc]{
\T{\H1\I3}\x\T{\J1\F2}\x\T{\H2\I4}\x\T{\J2}\\
\T{\L3}\T{\I2\H3}\T{\G3\L4}\T{\H4}\\
\T{}\x\T{\I2}\x\T{\L3}\T{\H3}\\
\T{\H2}\x\T{\J2}\T{}\T{}}
\\\vphantom{\dfrac{1}{2}}
x_1^3x_2x_3^2
\end{matrix}
$$
    \caption{the rest $12$ elements in $\HPD_{\tau}(\alpha)$}
    \label{fig:1302}
\end{figure}
\end{Eg}

\subsection{Branching rule}
\label{subsec:branch}

As an application of Theorem \ref{thm:mainthm}, let us give a two-side branching rule for key polynomials. 
Let $0\leq a<a+m\leq n$. 
We can give a combinatorial formula of the expansion 
$$\kappa_\alpha(x_1,\ldots,x_n)
=\sum_\beta c_{\alpha\beta}^{(a)}\cdot \kappa_{\beta}(x_{a+1},\ldots,x_{a+m})$$
where the sum over compositions $\beta\in \mathbb{Z}_{\geq 0}^{m}$ and 
$$c_{\alpha\beta}^{(a)}\in \mathbb{Z}[x_1,\ldots,x_{a},x_{a+m+1},\ldots,x_n].$$

For this purpose, we need slightly extend our model. 
Assume we are given two compositions 
$\alpha\in \mathbb{Z}_{\geq 0}^{n}$ and 
$\beta\in \mathbb{Z}_{\geq0}^m$. 
We pick an $N$ such that $\alpha_i\leq N$ and $\beta_i\leq N$ for all $i$.
We fix a choice of $\tau\in \{W,E\}^{n-m}$ with $a$ many $W$'s.
Consider a grid of $(n-m)\times (N+1)$ as in \eqref{eq:grid}. 
We can repeat what we did in section \ref{subsec:model}, except now we also label the color $i+a\in [n]$ on the $i$-th floor of the $\beta_{i}$-th column. 
We will denote the set of tilings by $\HPD_{\tau}(\alpha/\beta)$. 
For any $\pi\in \HPD_{\tau}(\alpha/\beta)$, we define its weight $\wt(\pi)$ by the same formula as in \eqref{eq:weight}. 

We have a simple bijection between
$$
\HPD_{\tau}(\alpha/\beta)\times 
\HPD_{\tau'}(\beta)\cong 
\HPD_{\tau\oplus \tau'}(\alpha)$$
by stacking two tilings. 
By Theorem \ref{thm:mainthm}, it follows immediately that 

\begin{Coro}\label{coro:branch}
For any choice of $\tau\in \{W,E\}^{n-m}$ with $a$ many $W$'s, 
we have 
$$c_{\alpha\beta}
=\sum_{\pi\in \HPD_{\tau}(\alpha/\beta)}\wt(\pi).$$
\end{Coro}

\begin{Eg}
Let us consider $\alpha=(1,3,0,2)$ and $\beta = (3,1)$ with $\tau=WE$. 
Here are the elements of $\HPD_{\tau}(\alpha/\beta)$: 
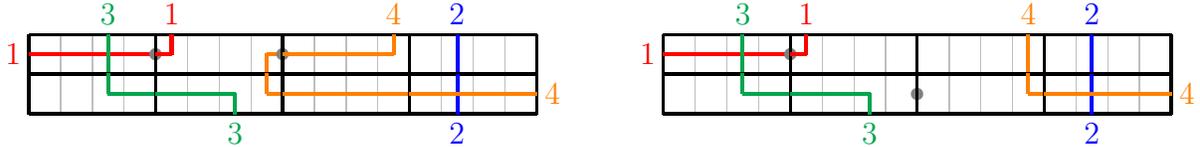
$$
\begin{matrix}
\PD[4]{
\M0{}\MM{{}{}3{}}\MM{1}\MM{{}{}{}4}\MM{{}2}\\
\M11\T{\H1\I3}\x\T{\J1\F4}\x\T{\J4}\T{\I2}\\
\M0{}\T{\L3}\T{\G3\L4}\T{\H4}\T{\I2\H4}\M44\\
\M0{}\MM{{}}\MM{{}{}3}\MM{{}}\MM{{}2}}
\end{matrix}\qquad 
\begin{matrix}
\PD[4]{
\M0{}\MM{{}{}3{}}\MM{1}\MM{{}{}{}4}\MM{{}2}\\
\M11\T{\H1\I3}\x\T{\J1}\T{\I4}\T{\I2}\\
\M0{}\T{\L3}\T{\G3}\x\T{\L4}\T{\I2\H4}\M44\\
\M0{}\MM{{}}\MM{{}{}3}\MM{{}}\MM{{}2}}
\end{matrix}
$$
We see that
$c_{\alpha\beta}=x_1^2+x_1x_4$. 
\end{Eg}

\subsection{Preparation for the proof}
\label{subsec:Preproof}

Let us discuss the proof of Theorem \ref{thm:mainthm}. 
When $\tau=W^n$, only the west tiles are used. 
By comparing our model with the weights in Borodin and Wheeler \cite[(3.7)]{BW22}, we see under the notation from loc. cit., 
the weighted sum 
$$\sum_{\pi\in \HPD_{\tau}(\alpha)}\wt(\pi)
= f^{n\cdots 21}_{\alpha}(x_1,\ldots,x_n;q=0,t=0).
$$
As pointed out under \cite[Theorem 1.3]{BW22}, the polynomial $f^{n\cdots 21}_\alpha$ is the permuted basement non-symmetric Macdonald polynomial introduced by Alexandersson \cite{Al19}. 
The specialization $q=t=0$ gives key polynomials, see the table under \cite[Question 7]{Al19}. 

Therefore, it reduces to finding a weight-preserving bijection among HPDs for different $\tau$'s. 
Following the idea of Knutson and Udell \cite{KU24}, it suffices to solve two basic cases:
\begin{enumerate}[\quad (a)]
    \item changing the type of the bottom row; \label{enm:Change}
    \item swapping two adjacent rows of different types.\label{enm:Swap}
\end{enumerate}
The case \eqref{enm:Change} is easy, simply by flipping the unique pipe in the bottom row:
$$
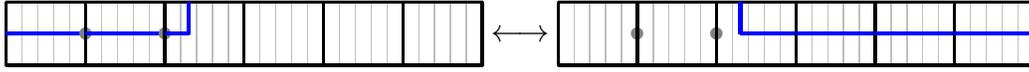
\PD[5][0.5pc][2pc]{
\T{\H2}\x\T{\H2}\x\T{\J2}\T{}\T{}\T{}
}
\longleftrightarrow 
\PD[5][0.5pc][2pc]{
\T{}\x\T{}\x\T{\L2}\T{\H2}\T{\H2}\T{\H2}
}
$$
The rest of paper will be devoted to construct a bijection for case \eqref{enm:Swap}. 

\begin{Eg}
Let us consider $n=3$ and $\alpha=(2,1,2)$. We choose $N=2$. 
Figure \ref{fig:212} illustrates \eqref{enm:Change} and \eqref{enm:Swap}.
\begin{figure}[hhhhhhhhh]
\[
\begin{matrix}
    \begin{matrix}
    \PD[1]{\M11\\\M22\\\M33}
        \PD[3][0.5pc]{
        \T{\H1}\x\T{\H1\I2}\x\T{\J1\I3}\\
        \T{\H2}\x\T{\J2}\T{\I3}\\
        \T{\H3}\x\T{\H3}\x\T{\J3}}\quad 
        \PD[3][0.5pc]{
        \T{\H1}\x\T{\H1\I2}\x\T{\J1\I3}\\
        \T{\H2}\x\T{\J2\F3}\x\T{\J3}\\
        \T{\H3}\x\T{\J3}\T{}}
    \PD[1]{\M0W\\\M0W\\\M0W}
    \end{matrix}
    &\leftarrow\eqref{enm:Change}\rightarrow&
    \begin{matrix}
    \PD[1]{\M11\\\M22\\\M0E}
        \PD[3][0.5pc]{
        \T{\H1}\x\T{\H1\I2}\x\T{\J1\I3}\\
        \T{\H2}\x\T{\J2}\T{\I3}\\
        \T{}\x\T{}\x\T{\L3}}\quad 
        \PD[3][0.5pc]{
        \T{\H1}\x\T{\H1\I2}\x\T{\J1\I3}\\
        \T{\H2}\x\T{\J2\F3}\x\T{\J3}\\
        \T{}\x\T{\L3}\T{\H3}}
    \PD[1]{\M0W\\\M0W\\\M33}
    \end{matrix}\\
    &&\uparrow\\
    &&\eqref{enm:Swap}\\
    &&\downarrow\\
    \begin{matrix}
    \PD[1]{\M11\\\M0E\\\M0E}
        \PD[3][0.5pc]{
        \T{\H1}\x\T{\H1\I2}\x\T{\J1\I3}\\
        \T{}\x\T{\I2}\x\T{\L3}\\
        \T{}\x\T{\L2}\T{\H2}}\quad
        \PD[3][0.5pc]{
        \T{\H1}\x\T{\H1\I2}\x\T{\J1\I3}\\
        \T{}\x\T{\L2}\T{\G2\L3}\\
        \T{}\x\T{}\x\T{\L2}}
    \PD[1]{\M0W\\\M33\\\M22}
    \end{matrix}
    &\leftarrow\eqref{enm:Change}\rightarrow & 
    \begin{matrix}
    \PD[1]{\M11\\\M0E\\\M22}
        \PD[3][0.5pc]{
        \T{\H1}\x\T{\H1\I2}\x\T{\J1\I3}\\
        \T{}\x\T{\I2}\x\T{\L3}\\
        \T{\H2}\x\T{\J2}\T{}}\quad
        \PD[3][0.5pc]{
        \T{\H1}\x\T{\H1\I2}\x\T{\J1\I3}\\
        \T{}\x\T{\L2}\T{\G2\L3}\\
        \T{\H2}\x\T{\H2}\x\T{\J2}}
    \PD[1]{\M0W\\\M33\\\M0W}
    \end{matrix}\\
    \uparrow&&\uparrow\\
    \eqref{enm:Swap}&&\eqref{enm:Swap}\\
    \downarrow&&\downarrow\\
    \begin{matrix}
    \PD[1]{\M0E\\\M11\\\M0E}
        \PD[3][0.5pc]{
        \T{}\x\T{\I2}\x\T{\I1\L3}\\
        \T{\H1}\x\T{\H1\I2}\x\T{\J1}\\
        \T{}\x\T{\L2}\T{\H2}}\quad
        \PD[3][0.5pc]{
        \T{}\x\T{\L2}\T{\I1\G2\L3}\\
        \T{\H1}\x\T{\H1}\x\T{\J1\I2}\\
        \T{}\x\T{}\x\T{\L2}}
    \PD[1]{\M33\\\M0W\\\M22}
    \end{matrix}
    &\leftarrow\eqref{enm:Change}\rightarrow & 
    \begin{matrix}
    \PD[1]{\M0E\\\M11\\\M22}
        \PD[3][0.5pc]{
        \T{}\x\T{\I2}\x\T{\I1\L3}\\
        \T{\H1}\x\T{\H1\I2}\x\T{\J1}\\
        \T{\H2}\x\T{\J2}\T{}}\quad
        \PD[3][0.5pc]{
        \T{}\x\T{\L2}\T{\I1\G2\L3}\\
        \T{\H1}\x\T{\H1}\x\T{\J1\I2}\\
        \T{\H2}\x\T{\H2}\x\T{\J2}}
    \PD[1]{\M33\\\M0W\\\M0W}
    \end{matrix}\\
    \uparrow\\
    \eqref{enm:Swap}\\
    \downarrow\\
    \begin{matrix}
    \PD[1]{\M0E\\\M0E\\\M11}
        \PD[3][0.5pc]{
        \T{}\x\T{\I2}\x\T{\I1\L3}\\
        \T{}\x\T{\L2}\T{\I1\H2}\\
        \T{\H1}\x\T{\H1}\x\T{\J1}}\quad
        \PD[3][0.5pc]{
        \T{}\x\T{\L2}\T{\I1\G2\L3}\\
        \T{}\x\T{}\x\T{\I1\L2}\\
        \T{\H1}\x\T{\H1}\x\T{\J1}}
    \PD[1]{\M33\\\M22\\\M0W}
    \end{matrix}
    &\leftarrow\eqref{enm:Change}\rightarrow & 
    \begin{matrix}
    \PD[1]{\M0E\\\M0E\\\M0E}
        \PD[3][0.5pc]{
        \T{}\x\T{\I2}\x\T{\I1\L3}\\
        \T{}\x\T{\L2}\T{\I1\H2}\\
        \T{}\x\T{}\x\T{\L1}}\quad
        \PD[3][0.5pc]{
        \T{}\x\T{\L2}\T{\I1\G2\L3}\\
        \T{}\x\T{}\x\T{\I1\L2}\\
        \T{}\x\T{}\x\T{\L1}}
    \PD[1]{\M11\\\M22\\\M33}
    \end{matrix}
\end{matrix}    
\]
    \caption{The cases \eqref{enm:Change} and \eqref{enm:Swap}}
    \label{fig:212}
\end{figure}
\end{Eg}

\section{Swapping two rows}\label{sec:swaprow}

In this section, we will establish the case \eqref{enm:Swap}. The bijection is constructed at the sub-tile level.
Namely, we will cut each column into $n$ sub-columns, resulting $(N+1)n$ sub-columns. 
Let $\pi \in \HPD_{\tau}(\alpha)$. 
Let us first state an equivalent way of computing weight $\wt(\pi)$. 
We define the weight of a west (resp., east) sub-tile with a pipe touching (resp., not touching) the right wall, i.e. 
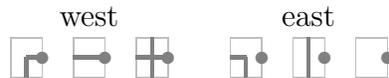
$$\begin{array}{c@{\qquad}c}
\text{west} & \text{east}\\
\PD{\F0\x}\quad 
\PD{\H0\x}\quad
\PD{\X00\x} & 
\PD{\G0\x}\quad
\PD{\I0\x}\quad
\PD{\O0{}\x}
\end{array}$$
to be $x_i$ if it is in the $i$-row and in the last sub-columns of a column.  
The weights of all other sub-tiles are $1$. 
Then the weight $\wt(\pi)$ is the product of the weights of sub-tiles in $\pi$.

Let us consider two adjacent rows of type $WE$
$$\PD[5]{
\M1{a}
\T{\O3\cdots}\T{\O3\cdots}
\MM{{}{}\cdots{}{}}
\T{\O3\cdots}\T{\O3\cdots}\\
\,
\T{\O3\cdots}\T{\O3\cdots}
\MM{{}{}\cdots{}{}}
\T{\O3\cdots}\T{\O3\cdots}\M4{b}
}$$
labeled by colors $a<b$ respectively. 
Note all the pipes appearing in the two rows are those pipes colored by $\{a,\ldots,b\}$. 
Namely, the lower boundary is labeled by colors $\{a+1,\ldots,b-1\}$ and the upper boundary is labeled by colors $\{a,\ldots,b\}$.
For a given boundary condition, we are going to establish a weight preserving bijection between type $WE$ and the  corresponding two rows of type $EW$: 
$$\PD[5]{
\,
\T{\O3\cdots}\T{\O3\cdots}
\MM{{}{}\cdots{}{}}
\T{\O3\cdots}\T{\O3\cdots}\M4{b}\\
\M1{a}
\T{\O3\cdots}\T{\O3\cdots}
\MM{{}{}\cdots{}{}}
\T{\O3\cdots}\T{\O3\cdots}
}$$
To do this, we need more notations.
For $a\leq c\leq b$, we define the source $s(c)$ (resp., the target $t(c)$) to be the sub-column colored by $c$ on the lower (resp., upper) boundary. 
Here we take the convention that define $s(a)=-\infty$ and $s(b)=\infty$. 

We say a sub-column is \emph{frozen} if it is between sub-columns $s(c)$ and $t(c)$ for some $a\leq c\leq b$ (not including $s(c)$ and $t(c)$).
Otherwise we say a sub-column is \emph{unfrozen}. 
We say a unfrozen sub-column is \emph{critical} if $s(c)$ or $t(c)$ for some $a\leq c\leq b$.

\begin{Eg}\label{eg:n=7}
Consider the following example when $n=7$. 
The following is a tiling of type $WE$. 
$$
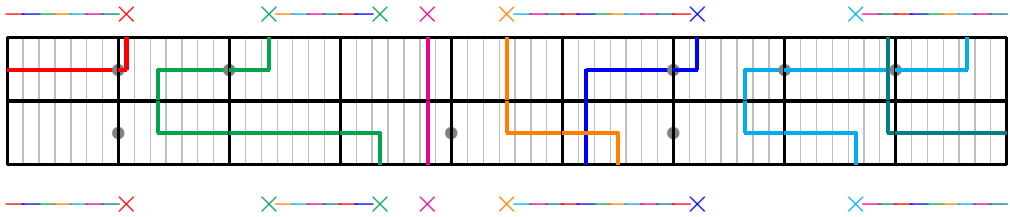
\PD[7][0.5pc][2pc]{
\MM{-------}\MM{\times}\MM{{}{}\times----}
\MM{--\times{}{}\times}
\MM{{}{}{}\times---}
\MM{-------}
\MM{-\times}\MM{{}{}{}{}\times--}\MM{-------}\\
\T{\H1}\x\T{\J1\F3}\x\T{\J3}\T{\I6}\T{\I4}\T{\F2}\x\T{\J2\F5}\x\T{\H5\I7}\x\T{\J5}\\
\T{}\x\T{\L3}\T{\H3}\T{\G3\I6}\x\T{\L4}\T{\I2\G4}\x\T{\L5}\T{\G5\L7}\T{\H7}\\
\MM{-------}\MM{\times}\MM{{}{}\times----}
\MM{--\times{}{}\times}
\MM{{}{}{}\times---}
\MM{-------}
\MM{-\times}\MM{{}{}{}{}\times--}\MM{-------}
}$$
We label the frozen sub-columns (by ``$-$'') and the critical sub-columns (by ``$\times$'') on it.
For readers who are eager to gain an insight into the bijection we will establish, we recommend referring directly to Example \ref{eg:n=7bi} later.
\end{Eg}

\subsection{Frozen sub-columns}
\label{subsec:frozen}
Let us first deal with a single frozen sub-column. 
Assume it is between $s(c)$ and $t(c)$ for $a\leq c\leq b$. 
Since in rows of different types the pipe goes along different directions, the pipe colored by $c$ must appear in the frozen sub-column and can only go horizontally. 

For type $WE$, a frozen sub-column must be one of: 
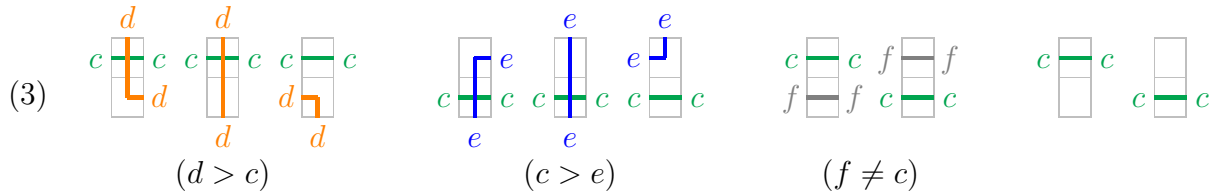
\begin{equation}\label{eq:frozentiles}
\begin{array}{c@{\qquad}c@{\qquad}c@{\qquad}c}
\PD[1]{\,\M4d\\    \M3c\X34\M3c\\    \,\L4\M4d\\    \,}
\PD[1]{\,\M4d\\    \M3c\X34\M3c\\    \,\I4\\    \,\M4d}
\PD[1]{\,\\    \M3c\H3\M3c\\    \M4d\G4\\    \,\M4d}&
\PD[1]{\,\\    \,\F2\M2e\\    \M3c\X32\M3c\\    \,\M2e}
\PD[1]{\,\M2e\\    \,\I2\\    \M3c\X32\M3c\\    \,\M2e}
\PD[1]{\,\M2e\\    \M2e\J2\\    \M3c\H3\M3c\\    \,}&
\PD[1]{\,\\\M3c\H3\M3c\\\M0f\H0\M0f\\    \,}
\PD[1]{\,\\\M0f\H0\M0f\\\M3c\H3\M3c\\    \,}&
\PD[1]{\,\\\M3c\H3\M3c\\\,\O0{}\,\\    \,}
\PD[1]{\,\\    \,\O0{}\\    \M3c\H3\M3c\\    \,}\\
(d>c)&(c>e)&(f\neq c)&
\end{array}
\end{equation}
Let us explain the condition of $d>c$, the condition for $c>e$ is similar. 
In the first and the second cases, we have to assume $d>c$ by our assumption on $\PD[1]{\X00}$. 
In the third case, the pipe colored by $c$ and $d$ will finally intersect, i.e. we will meet 
$$
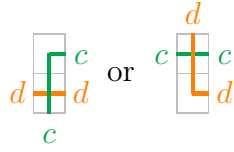
\PD[1]{\,\\    \,\F3\M3c\\    \M4d\X43\M4d\\    \,\M3c}
\text{ or }
\PD[1]{\,\M4d\\    \M3c\X34\M3c\\    \,\L4\M4d\\    \,}
$$
on its left. 
So we also have $d>c$.

For type $EW$, similar argument shows the sub-column must be one of
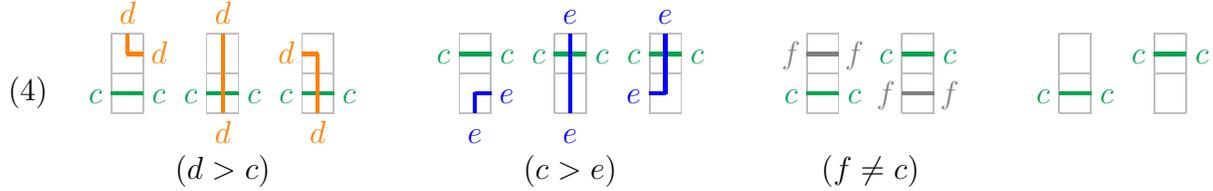
\begin{equation}\label{eq:frozentilesEW}
\begin{array}{c@{\qquad}c@{\qquad}c@{\qquad}c}
\PD[1]{\,\M4d\\    \,\L4\M4d\\    \M3c\H3\M3c\\    \,}
\PD[1]{\,\M4d\\    \,\I4\\    \M3c\X34\M3c\\    \,\M4d}
\PD[1]{\,\\    \M4d\G4\\    \M3c\X34\M3c\\    \,\M4d}&
\PD[1]{\,\\    \M3c\H3\M3c\\    \,\F2\M2e\\    \,\M2e}
\PD[1]{\,\M2e\\    \M3c\X32\M3c\\    \,\I2\\    \,\M2e}
\PD[1]{\,\M2e\\    \M3c\X32\M3c\\    \M2e\J2\\    \,}&
\PD[1]{\,\\\M0f\H0\M0f\\\M3c\H3\M3c\\\,}
\PD[1]{\,\\    \M3c\H3\M3c\\    \M0f\H0\M0f\\ \,}&
\PD[1]{\,\\    \,\O0{}\,\\    \M3c\H3\M3c\\    \,}
\PD[1]{\,\\    \M3c\H3\M3c\\    \,\O0{}\\    \,}\\
(d>c)&(c>e)&(f\neq c)&
\end{array}
\end{equation}
In a frozen sub-column, we define a bijection 
by sending sub-columns in \eqref{eq:frozentiles} to those in \eqref{eq:frozentilesEW} respectively. 
Note that colors on the left and right boundaries are swapped. 
It is straightforward to see this is a weight-preserving bijection. 

\subsection{Unfrozen sub-columns}
\label{subsec:unfrozen}
We still need to deal with unfrozen sub-columns. 
We say a critical sub-column is of type $M$ (resp., $L$, $R$) if both of its left and its right (resp., only its left, only its right) are unfrozen. 

Let us consider a maximal interval $I$ of unfrozen sub-columns, i.e. it is the collection of unfrozen sub-columns between two critical sub-columns of type $L$ and $R$ (including both critical sub-columns).
The critical sub-columns $m_0,m_1,\ldots,m_k$ (from left to right) in $I$ cut the interval into sub-intervals $I_1,\ldots,I_k$ (from left to right).
We show them in a diagram
$$
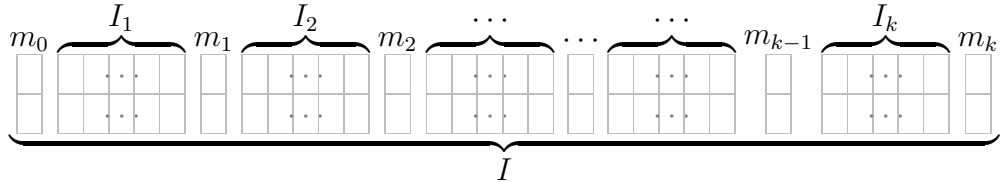
\def\c#1{%
\stackrel{\displaystyle #1}{\PD[1][0.8pc]{\O0{}\\\O0{}}}}
\def\I#1{%
\overbrace{\PD[1][0.8pc]{
\O0{}\O0{}\O0{\cdots}\O0{}\O0{}\\
\O0{}\O0{}\O0{\cdots}\O0{}\O0{}
}}^{\displaystyle #1}}
\underbrace{
\c{m_0}\I{I_1}
\c{m_1}\I{I_2}
\c{m_2}\I{\cdots}
\c{\cdots}\I{\cdots}
\c{m_{k-1}}\I{I_k}\c{m_k}
}_{\displaystyle I}
$$
Let $N_i$ ($i=1,\ldots,k$) be the maximal integer such that 
$nN_i\leq |I_i|$. 
We will construct a bijection between
tilings of the interval $I$ of type $WE$ (resp., $EW$)  and the set 
\begin{equation}\label{eq:unfrozenbijto}
\bigg\{(j_1,\ldots,j_k): 
0\leq j_i\leq N_i \text{ for $i=1,\cdots,k$}
\bigg\}.
\end{equation}

Let us first deal with the type $WE$. 
Let $1\leq i\leq k$. 
In the sub-columns $I_i\cup\{m_i\}$, there must exists a sub-column $m$ such that the tiling is 
$$
\PD[5]{
\M3{}\M3{}\MM{{}}
\MM{{}{\color{black}m}}\M0{}\M2{\cdots}\M0{}\MM{{}}\M0{}\M2{\color{black}m_i}\\
\O3{}\O3{}\T{}\T{\F2}\M0{}\M2{\cdots}\M0{}\T{\H2}\H2\J2\\
\O3{}\O3{}\T{}\T{\L2}\M0{}\M2{\cdots}\M0{}\T{\H2}\H2\H2}
\quad \text{ or }
\PD[5]{
\M0{}\M0{}\M2{\color{black}m_i}\\
\M0{\cdots}\M0{}\H2\\
\M0{\cdots}\M0{}\G2}
\quad \text{ or }
\PD[5]{
\M0{}\M0{}\M2{\color{black}m_i}\\
\M0{\cdots}\M0{}\J2\\
\M0{\cdots}\M0{}\G2}
$$
depending on the upper and lower boundary conditions on the $m_i$-th sub-column.
When $m_i=m$, the tiling is understood as 
$$\PD[5]{
\M3{}\M3{}\MM{{}}
\MM{{}{}}\M0{}\M2{\cdots}\M0{}\MM{{}}\M0{}\M2{\color{black}m_i}\\
\O3{}\O3{}\T{}\T{}\M0{}\M2{\cdots}\M0{}\T{}\O0{}\I2\\
\O3{}\O3{}\T{}\T{}\M0{}\M2{\cdots}\M0{}\T{}\O0{}\L2}
\quad \text{ or }
\PD[5]{
\M0{}\M0{}\M2{\color{black}m_i}\\
\M2{\cdots}\M0{}\F2\\
\M2{\cdots}\M0{}\I2}
\quad \text{ or }
\PD[5]{
\M0{}\M0{}\M2{\color{black}m_i}\\
\M2{\cdots}\M0{}\I2\\
\M2{\cdots}\M0{}\I2}
$$
Let us define $j_i=\dfrac{m_i-m}{n}$. 
Since $m\in I_i\cup \{m_i\}$, we have $0\leq j_i\leq N_i$. 
This defines a map to \eqref{eq:unfrozenbijto}. 
Note that each $j_i$ determines the tiling in the interval $I_i\cup \{m_i\}$. Moreover on $m_0$, the tiling must be 
$$\PD[1]{\J2\\\O0{}}
\quad\text{ or }\quad
\PD[1]{\O0{}\\\G2}$$
depending on the upper and lower boundary conditions on the $m_0$-th sub-column.
This proves the map is a bijection. 

Similar consideration applies to the type $EW$. 
In the sub-columns $\{m_{i-1}\}\cup i_{i}$, there exists a sub-column $m$ such that the tiling is 
$$
\PD[5]{
\M4{\color{black}m_{i-1}}\M0{}\MM{{}}\M0{}\M4{\cdots}\M0{}\MM{{}{}{}{\color{black}m}}\\
\L4\H4\T{\H{4}}\M0{}\M4{\cdots}\M0{}\T{\G{4}}\O0{}\O0{}\\
\H4\H4\T{\H{4}}\M0{}\M4{\cdots}\M0{}\T{\J{4}}\O0{}\O0{}\\
}
\quad \text{ or }\quad 
\PD[5]{
\M4{\color{black}m_{i-1}}\M0{}\\
\H4\M0{}\M4{\cdots}\\
\F4\M0{}\M4{\cdots}}
\text{ or }\quad 
\PD[5]{
\M4{\color{black}m_{i-1}}\M0{}\\
\L4\M0{}\M4{\cdots}\\
\F4\M0{}\M4{\cdots}}
$$
depending the boundary condition on the $m_{i-1}$-th sub-column.
When $m_{i-1}=m$, the tiling is understood as 
$$
\PD[5]{
\M4{\color{black}m_{i-1}}\M0{}\MM{{}}\M0{}\M4{\cdots}\M0{}\MM{{}{}{}{}}\\
\I4\O0{}\T{}\M0{}\M4{\cdots}\M0{}\T{}\O0{}\O0{}\\
\J4\O0{}\T{}\M0{}\M4{\cdots}\M0{}\T{}\O0{}\O0{}\\
}\quad
\text{ or }\quad 
\PD[5]{
\M4{\color{black}m_{i-1}}\M0{}\\
\G4\M0{}\M4{\cdots}\\
\I4\M0{}\M4{\cdots}}
\text{ or }\quad 
\PD[5]{
\M4{\color{black}m_{i-1}}\M0{}\\
\I4\M0{}\M4{\cdots}\\
\I4\M0{}\M4{\cdots}}$$
We define $j_i=\dfrac{m-m_{i-1}}{n}$.
Similarly as the $WE$ cases, this defines a bijection to \eqref{eq:unfrozenbijto}. 

\begin{Eg}
Let us consider the following tiling of type $WE$
$$\PD[5][0.7pc][2pc]{
\M3{\times}\M0{}\M0{}
\MM{{}}\MM{{}{}{}\times}\MM{{}}\MM{{}}\MM{{}{}{}{}\times}
\MM{{}}\MM{{}\times}\\
\J3\O0{}\O0{}
\T{\F4}\x\T{\J4}\T{\F5}\x\T{\H5}\x\T{\J5}\T{\F2}\x\H2\J2\\
\O0{}\O0{}\O0{}
\x\T{\L4}\T{\G4}\x\T{\L5}\T{\H5}\T{\G5}\x\T{\L2}\H2\H2\\
\M3{\times}\M0{}\M0{}
\MM{{}}\MM{{}{}{}\times}\MM{{}}\MM{{}}\MM{{}{}{}{}\times}
\MM{{}}\MM{{}\times}}$$
It corresponds to $(1,2,1)$. 
The corresponding tiling of type $EW$ is 
$$\PD[5][0.7pc][2pc]{
\M3{\times}\M0{}\M0{}
\MM{{}}\MM{{}{}{}\times}\MM{{}}\MM{{}}\MM{{}{}{}{}\times}
\MM{{}}\MM{{}\times}\\
\L3\H3{}\H3{}
\T{\G3}\x\T{\L4}\T{\H4}\T{\G4}\x\T{\L5}\T{\G5}\x\O2{}\L2\\
\H3\H3\H3\x
\T{\J3}\T{\F4}\x\T{\H4}\x\T{\J4}\T{\F5}\x\T{\J5}\O2{}\O2{}\\
\M3{\times}\M0{}\M0{}
\MM{{}}\MM{{}{}{}\times}\MM{{}}\MM{{}}\MM{{}{}{}{}\times}
\MM{{}}\MM{{}\times}}$$
\end{Eg}

Now we get a bijection 
between the tilings of the interval $I$ of type $WE$ and of type $EW$. 
From the discussion above, it is direct to see the colors on the boundary of $I$ are swapped. 
Let $w$ be the number of sub-columns in $I\setminus \{m_k\}$ appearing as the last sub-columns in a column.
Then it is clear that the weight of a tiling is 
$$x_a^{j_1+\cdots+j_k}x_b^{w-j_1-\cdots-j_k}$$
for both $WE$ and $EW$. 
This shows the bijection is weight-preserving. 

\begin{Eg}\label{eg:n=7bi}
Let us apply our bijection to Example \ref{eg:n=7}. 
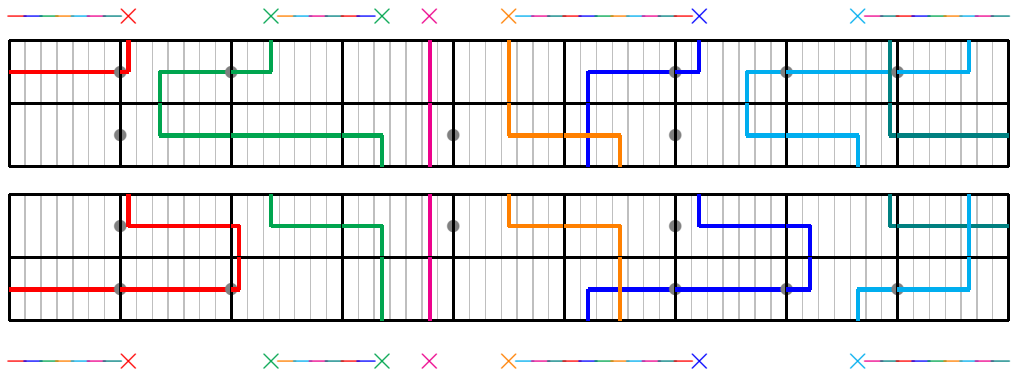
$$
\PD[7][0.5pc][2pc]{
\MM{-------}\MM{\times}\MM{{}{}\times----}
\MM{--\times{}{}\times}
\MM{{}{}{}\times---}
\MM{-------}
\MM{-\times}\MM{{}{}{}{}\times--}\MM{-------}\\
\T{\H1}\x\T{\J1\F3}\x\T{\J3}\T{\I6}\T{\I4}\T{\F2}\x\T{\J2\F5}\x\T{\H5\I7}\x\T{\J5}\\
\T{}\x\T{\L3}\T{\H3}\T{\G3\I6}\x\T{\L4}\T{\I2\G4}\x\T{\L5}\T{\G5\L7}\T{\H7}\\[2ex]\\
\T{}\x\T{\L1}\T{\G1\L3}\T{\G3\I6}\x\T{\L4}\T{\G4}\x\T{\L2}\T{\G2\L7}\T{\H{7}\I5}\\
\T{\H1}\x\T{\H1}\x\T{\J1}\T{\I3\I6}\T{}\T{\F2\I4}\x\T{\H2}\x\T{\J2\F5}\x\T{\J5}\\
\MM{-------}\MM{\times}\MM{{}{}\times----}
\MM{--\times{}{}\times}
\MM{{}{}{}\times---}
\MM{-------}
\MM{-\times}\MM{{}{}{}{}\times--}\MM{-------}
}$$
\end{Eg}

Combining the two bijections for frozen and unfrozen sub-columns together, we get a bijection of two adjacent rows of type $WE$ and type $EW$. 
The case \eqref{enm:Swap} is established.

\end{document}